\documentclass[12pt,twoside,english]{svjour3}
\usepackage[T1]{fontenc}
\usepackage[latin9]{inputenc}
\usepackage[a4paper]{geometry}
\usepackage{float}
\usepackage{url}
\usepackage{enumitem}
\usepackage{amsmath}
\usepackage{amssymb}
\usepackage{graphicx}
\usepackage{setspace}
\PassOptionsToPackage{normalem}{ulem}
\usepackage{ulem}
\makeatletter

\floatstyle{ruled}
\newfloat{algorithm}{tbp}{loa}
\providecommand{\algorithmname}{Algorithm}
\floatname{algorithm}{\protect\algorithmname}

  \newenvironment{svmultproof}{\begin{proof}}{\qed\end{proof}}

\@ifundefined{date}{}{\date{}}
\makeatother

\usepackage{babel}

\begin{document}
\begin{onehalfspace}

\title{Convergence Analysis of Processes with Valiant Projection Operators
in Hilbert Space}
\end{onehalfspace}

\author{Yair Censor\textsuperscript{1} and Rafiq Mansour\textsuperscript{2}}

\authorrunning{Y. Censor and R. Mansour}

\institute{\textsuperscript{1}Corresponding Author. Department of Mathematics,
University of Haifa, Mt. Carmel, Haifa 3498838, Israel. E-mail: yair@math.haifa.ac.il.
\textsuperscript{2}Department of Mathematics, University of Haifa,
Mt. Carmel, Haifa 3498838, Israel. E-mail: intogral@gmail.com.}
\maketitle
\begin{center}
December 4, 2016. Revised: May 29, 2017. Revised: September 11, 2017.
\par\end{center}
\begin{abstract}
Convex feasibility problems require to find a point in the intersection
of a finite family of convex sets. We propose to solve such problems
by performing set-enlargements and applying a new kind of projection
operators called valiant projectors. A valiant projector onto a convex
set implements a special relaxation strategy, proposed by Goffin in
1971, that dictates the move toward the projection according to the
distance from the set. Contrary to past realizations of this strategy,
our valiant projection operator implements the strategy in a continuous
fashion. We study properties of valiant projectors and prove convergence
of our new valiant projections method. These results include as a
special case and extend the 1985 automatic relaxation method of Censor.
\end{abstract}

\keywords{Intrepid projector; convex feasibility problem; valiant projector; set-enlargement; automatic relaxation method (ARM); ART3 algorithm; Goffin's principle} \subclass{65K05, 90C25}

\numberwithin{definition}{section}

\numberwithin{proposition}{section}

\numberwithin{remark}{section}

\numberwithin{theorem}{section}

\numberwithin{lemma}{section}

\numberwithin{corollary}{section}

\numberwithin{algorithm}{section}

\renewcommand{\labelenumi}{\roman{enumi}.}

\section{Introduction}

We consider the convex feasibility problem (CFP) in a real Hilbert
space. It consists of a finite family of closed and convex sets with
a nonempty intersection and calls to find an element in this intersection.
There are many algorithms in the literature for solving CFPs, see,
e.g., \cite{bb96,censor2015projection}, and many problems in operations
research and in various other fields can be presented as feasibility
problems. For example, a problem of road design is modeled as a feasibility
problem in \cite{bauschke2015projection}, where the motivation is
automated design of road alignments. A road alignment is represented
by the center-line of the road, which is idealized as a (generally)
nonlinear, smooth curve. To facilitate construction drawings, civil
engineers reduce the three-dimensional road design to two two-dimensional
parts, horizontal and vertical. For a new approach to road design
optimization see \cite{bauschke2015stadium}. In \cite{artacho2013douglas}
the authors give general recommendations for successful application
of the Douglas\textendash Rachford feasibility-seeking method to convex
and non-convex real matrix-completion problems. The work in \cite{borwein2014reflection}
focuses on the problem of protein conformation determination formulated
within the framework of matrix completion by solving CFPs. In \cite[section 6]{artacho2013recent}
solving Sudoku puzzles is modeled as an integer feasibility problem.
The books \cite{censor1997parallel}, \cite{Ceg-book} and \cite{BC11}
contain many algorithms and methods that solve the CFP, and there
is also a description of a wide range of operators and their properties.
Hence, the importance of the CFP stems from its flexibility to accommodate
problems from various fields and simplify their understanding and
solutions. 

Searching for a solution to a system of linear equations is a convex
feasibility problem and has led to many different iterative methods.
When the system of linear equations is inconsistent, due to modeling
or measurements inaccuracies, it has been suggested to replace it
by a system of pairs of opposing linear inequalities, that represent
nonempty hyperslabs. Applying projection methods to this problem can
be done by using any iterative method for linear inequalities, such
as the method of Agmon \cite{agmon1954relaxation} and Motzkin and
Schoenberg \cite{motzkin1954relaxation} (AMS). However, in order
to improve computational efficiency, Goffin \cite{goffin1971finite}
proposed to replace projections onto the hyperslabs by a strategy
of projecting onto the original hyperplane (from which the hyperslab
was created), when the current iterate is ``far away'' from the
hyperslab, and reflecting into the hyperslab's boundary, when the
current iterate is ``close to the hyperslab'', while keeping the
iterate unchanged if it is already inside the hyperslab.

In \cite{herman1975relaxation} Herman suggested to implement Goffin's
strategy by using an additional enveloping hyperslab in order to determine
the ``far'' and the ``close'' distance of points from the hyperplane,
resulting in his ``Algebraic Reconstruction Technique 3'' (ART3)
algorithm. In \cite{censor1985automatic} Censor also embraced the
idea of hyperslabs, and defined an algorithmic operator that implemented
Goffin's strategy in a continuous manner, resulting in the Automatic
Relaxation Method (ARM). For applications and additional details see
\cite{censor1987some} and \cite{censor1988parallel}. 

A fundamental question, that remained open since then, was whether
the hyperslabs approach to handle linear equations and Goffin's principle
can be applied to general convex sets and not only to linear equations.
This question was recently studied by Bauschke, Iorio and Koch in
\cite{bauschke2014method}, see also \cite{bauschke2015projection}
and \cite{bauschke2015stadium} for further details and interesting
applications. They defined convex sets enlargements instead of hyperslabs
and used them to generalize the algorithmic operator that appeared
in \cite{herman1975relaxation}. They defined an operator which they
called the intrepid projector, intended to generalize the ART3 algorithm
of \cite{herman1975relaxation} to convex sets. Motivated by \cite{bauschke2014method},
we present in this paper a new operator which we call the valiant
operator, that enables to implement the algorithmic principle embodied
in the ARM of \cite{censor1985automatic} to general convex feasibility
problems. Observe that both ART3 and ARM seek a feasible point in
the intersections of the hyperslabs and so their generalizations to
the convex case seek feasibility of appropriate enlargement sets that
define the extended problem.

The new valiant projection method (VPM) proposed and studied here
answers affirmatively the theoretical question posed above. To date
we have no computational experience with it, that will allow us make
any claims about its actual advantages. It may be the case that using
valiant operators is beneficial not always but only in some specific
situations (size, sparsity, nature of problems, parameters, specific
applications, etc.) or under some additional conditions. To discover
these a methodological numerical work is required and we plan to undertake
such work with collaborators from some application fields with which
we are involved.

The paper is organized as follows: In Section \ref{sec:Preliminaries},
we give definitions and preliminaries. In Section \ref{sec:Motivation},
we present the motivation of the main idea of this paper and in Section
\ref{sec:The-valiant-projector}, we present the new algorithmic projector
and its relevant features. Finally, in Section \ref{sec:The-method},
we present our new algorithm and prove its convergence, and in Section
\ref{sec:Conclusions} we offer concluding comments.

\section{Preliminaries\label{sec:Preliminaries}}

For the reader's convenience we include in this section some properties
of operators in Hilbert space that will be used to prove our results.
We use the recent excellent book of Cegielski \cite{Ceg-book} as
our desk-copy in which all the results of this section can be found
\cite[Chapter 2 and Chapter 3]{Ceg-book}. Let $\mathcal{H}$ be a
real Hilbert space with inner product $\left\langle \cdot,\cdot\right\rangle $
and induced norm $\parallel\cdot\parallel$, and let $X\subseteq\mathcal{H}$
be a closed convex subset. Denote the index set $I:=\left\{ 1,2,\ldots,m\right\} $.
If $\varOmega\subseteq\mathcal{H}$ and $x\in\mathcal{H}$ then we
denote by $P_{\varOmega}(x)$ the metric projection of $x$ onto $\varOmega$.
\begin{definition}
\label{def:FM+NE+FNE}An operator $T:X\rightarrow\mathcal{H}$ is:
\begin{enumerate}
\item \textit{\emph{Nonexpansive}}\emph{ }(NE)\emph{,} if $\Vert T(x)-T(y)\Vert\leq\Vert x-y\Vert$
for all $x,y\in X$.
\item \textit{\emph{Firmly nonexpansive}} (FNE)\emph{,} if $\Vert T(x)-T(y)\Vert^{2}+\Vert(x-T(x))-(y-T(y))\Vert^{2}\leq\Vert x-y\Vert^{2}$
for all $x,y\in X.$
\end{enumerate}
\end{definition}
\begin{remark}
\label{rem:FNE=00003DNE}(i) It is clear from Definition \ref{def:FM+NE+FNE}
that every FNE operator is NE. See also, \cite[Theorem 2.2.10 (i)-(ii)]{Ceg-book}.
(ii) By \cite[Lemma 2.1.12]{Ceg-book} the family of NE is closed
under convex combinations and compositions.
\end{remark}
\begin{definition}
\label{def:projection onto a set} Let $C$ be a nonempty closed convex
subset of $\mathcal{H}$, let $x\in\mathcal{H}$, and let $c\in C$.
Denote the distance from $x$ to $C$ by $d_{C}(x):=\textup{inf}_{c\in C}\Vert x-c\Vert$,
the infimum is attained at a unique vector called the projection of
$x$ onto $C$ and denoted by $P_{C}(x)$.
\end{definition}
\begin{proposition}
\label{prop:projection=00003DFNE}\textup{\cite[Proposition 4.8]{BC11}}
Let $B$ be a nonempty closed convex subset of $\mathcal{H}$. Then
the projector $P_{B}$ is FNE.
\end{proposition}
\begin{definition}
\label{def:Fixed point-operator}An operator $T:X\rightarrow\mathcal{H}$
having a nonempty fixed point set $\textup{Fix}T=\left\{ x\in X\mid T(x)=x\right\} $
is:

\begin{enumerate}
\item \textit{\emph{Quasi-nonexpansive}} (QNE) if $\Vert T(x)-z\Vert\leq\Vert x-z\Vert$
for all $x\in X$ and $z\in$Fix$T$.
\item \textit{\emph{Strictly quasi-nonexpansive}} (sQNE) if $\Vert T(x)-z\Vert<\Vert x-z\Vert$
for all $x\notin$Fix$T$ and $z\in$Fix$T$.
\item \emph{B}\textit{\emph{-strictly quasi-nonexpansive}} (\textit{B}-sQNE),
where $B\neq\textrm{Ø}$ and $B\subseteq$Fix$T$, if $T$ is quasi-nonexpansive
and $\Vert T(x)-z\Vert<\Vert x-z\Vert$ for all $x\notin$Fix$T$
and $z\in B$.
\item \emph{$\alpha$-}\textit{\emph{strongly quasi-nonexpansive}} ($\alpha$-SQNE)
if $\Vert T(x)-z\Vert^{2}\leq\Vert x-z\Vert^{2}-\alpha\Vert T(x)-x\Vert^{2}$
for all $x\in X$ and $z\in$Fix$T$, where $\alpha\geq0$. If $\alpha>0$
then $T$ is called \textit{\emph{strongly quasi-nonexpansive}} (SQNE).
\end{enumerate}
\end{definition}
The next implications follow directly from the definitions, see \cite[page 47]{Ceg-book}
and \cite[Remark 2.1.44(iii)]{Ceg-book}.
\begin{proposition}
\label{prop:sQNE}For an operator $T:X\rightarrow\mathcal{H}$ having
a fixed point, the following statements hold: 

\begin{enumerate}
\item If $T$ is sQNE then $T$ is $B$-sQNE, where $B\subseteq$\textup{Fix}$T$. 
\item If $T$ is \textup{Fix}$T$-sQNE then $T$ is sQNE.
\item If $T$ is SQNE then it is sQNE.
\end{enumerate}
\end{proposition}
The following proposition presents the relationship between NE and
QNE operators.
\begin{proposition}
\label{prop:NE to QNE}\textup{\cite[Lemma 2.1.20]{Ceg-book}} An
NE operator $U:X\rightarrow\mathcal{H}$ with a fixed point is QNE.
\end{proposition}
\begin{remark}
\label{rem:SQNE=00003Dcomp+comb}From \cite[Corollary 2.1.47]{Ceg-book}
and \cite[Fig. 2.14]{Ceg-book}, a family of SQNE operators with a
common fixed point is closed under convex combinations and compositions.
\end{remark}
The following proposition shows that the relaxation of a projection
onto a nonempty closed convex set is an SQNE operator.
\begin{proposition}
\textup{\cite[Fact 1]{bauschke2014method}}\label{prop:relaxed projector}
Let $\varOmega$ be nonempty closed convex subset of $\mathcal{H}$
and let $\lambda\in]0,2[$. Set $R:=(1-\lambda)\textup{Id}+\lambda P_{\varOmega}$,
let $x\in\mathcal{H}$ and $c\in\varOmega$. Then 
\begin{equation}
\Vert x-c\Vert^{2}-\Vert R(x)-c\Vert^{2}\geq\frac{2-\lambda}{\lambda}\Vert x-R(x)\Vert^{2}.
\end{equation}
\end{proposition}
\begin{theorem}
\label{thm:Fix=00003DintersectionFix}\textup{\cite[Theorem 2.1.26(ii)]{Ceg-book}}
Let the operators $U_{i}:X\rightarrow X$ , $i\in I$, with $\bigcap_{i\in I}$\textup{Fix}$U_{i}\neq\textrm{Ø}$,
be $B$-sQNE, where $B\subseteq\bigcap_{i\in I}$\textup{Fix}$U_{i}$,
$B\neq\textrm{Ø}$. If $U:=U_{m}U_{m-1}\ldots U_{1}$ then
\end{theorem}
\textit{
\begin{equation}
\textnormal{Fix}U=\bigcap_{i\in I}\textnormal{Fix}U_{i},
\end{equation}
}

\noindent \textit{and $U$ is $B$-sQNE.}
\begin{definition}
\label{def:A-R-operator}An operator $U:X\rightarrow X$ is asymptotically
regular if for all $x\in X$,
\end{definition}
\noindent 
\begin{equation}
\underset{k\rightarrow\infty}{\textnormal{lim}}\Vert U^{k+1}(x)-U^{k}(x)\Vert=0.
\end{equation}

\begin{theorem}
\label{thm:SQNE=00003DAR}\textup{\cite[Theorem 3.4.3]{Ceg-book}}
Let $U:X\rightarrow X$ be an operator with a fixed point. If $U$
is SQNE then it is asymptotically regular.
\end{theorem}
The following well-known theorem is due to Opial.
\begin{theorem}
\label{thm:Opial}\textup{\cite[Theorem 3.5.1]{Ceg-book}} Let $X\subseteq\mathcal{H}$
be a nonempty closed convex subset of a Hilbert space $\mathcal{H}$
and let $U:X\rightarrow X$ be a nonexpansive and asymptotically regular
operator with a fixed point. Then, for any $x\in X$, the sequence
$\left\{ U^{k}(x)\right\} _{k=0}^{\infty}$ converges weakly to a
point $z\in\textup{Fix}U$.
\end{theorem}
\begin{definition}
\label{def:demi-closed}An operator $T:X\rightarrow\mathcal{H}$ is
demiclosed at 0 if for any weakly convergent sequence $x^{k}\rightharpoonup y\in X$
with $T(x^{k})\rightarrow0$ we have $T(y)=0$.
\end{definition}
The next theorem is known as the demiclosedness principle.
\begin{theorem}
\label{thm:demiclosedness principle}\textup{\cite[Lemma 3.2.5]{Ceg-book}
}Let $T:X\rightarrow\mathcal{H}$ be an NE operator and $y\in X$
be a weak cluster point of a sequence $\left\{ x^{k}\right\} _{k=0}^{\infty}$.
If $\left\Vert T\left(x^{k}\right)-x^{k}\right\Vert \longrightarrow0$,
then $y\in\textup{Fix}T$.
\end{theorem}
The following definition extends Definition \ref{def:A-R-operator}
to a sequence of operators.
\begin{definition}
\label{def:asymp-reg-seq}Let $X\subseteq\mathcal{H}$ be a nonempty
closed convex subset. We say that a sequence of operators $U_{k}:X\rightarrow X$
is asymptotically regular, if for any $x\in X$ 
\begin{equation}
\lim_{k}\left\Vert U_{k}U_{k-1}\ldots U_{0}\left(x\right)-U_{k-1}\ldots U_{0}\left(x\right)\right\Vert =0,
\end{equation}
or, equivalently, 
\begin{equation}
\lim_{k}\left\Vert U_{k}\left(x^{k}\right)-x^{k}\right\Vert =0,
\end{equation}
where the sequence $\left\{ x^{k}\right\} _{k=0}^{\infty}$ is generated
by recurrence $x^{k+1}=U_{k}\left(x^{k}\right)$ with $x^{0}=x$.
\end{definition}
\begin{theorem}
\label{thm:3.6.2}\textup{\cite[Theorem 3.6.2(i)]{Ceg-book}} Let
$X\subseteq\mathcal{H}$ be a nonempty closed convex subset, $S:X\rightarrow\mathcal{H}$
be an operator with a fixed point and such that $S-\textup{Id}$ is
demiclosed at 0. Let $\left\{ U_{k}\right\} _{k=0}^{\infty}$ be an
asymptotically regular sequence of quasi-nonexpansive operators $U_{k}:X\rightarrow X$
such that $\bigcap_{k=0}^{\infty}\textnormal{Fix}U_{k}\supseteq\textnormal{Fix}S$.
Let the sequence $\left\{ x^{k}\right\} _{k=0}^{\infty}$ be generated
by recurrence $x^{k+1}=U_{k}\left(x^{k}\right)$, with an arbitrary
$x^{0}\in X$. If the sequence of operators $\left\{ U_{k}\right\} _{k=0}^{\infty}$
has the property 
\begin{equation}
\lim_{k}\left\Vert U_{k}\left(x^{k}\right)-x^{k}\right\Vert =0\;\;\;\;\Longrightarrow\;\;\;\;\lim_{k}\left\Vert S\left(x^{k}\right)-x^{k}\right\Vert =0,\label{eq:6}
\end{equation}
then $\left\{ x^{k}\right\} _{k=0}^{\infty}$ converges weakly to
a point $z^{*}\in\textnormal{Fix}S$.
\end{theorem}

\section{The Valiant Projector: Intuition and Motivation\label{sec:Motivation}}

In this section, we explain the motivation and intuition behind valiant
projection operators which relies on two basic ideas: the notion of
enlargement of a convex set, and a strategy proposed by Goffin \cite{goffin1971finite}
that dictates to move towards the projection according to the distance
to the set. In more details, the development was as follows. First,
Herman proposed and studied in \cite{herman1975relaxation} the ``Algebraic
Reconstruction Technique 3'' (ART3) algorithm for solving a system
of two-sided linear inequalities. He set out to solve a large and
sparse, possibly inconsistent, system of linear equations stemming
from the problem of image reconstruction from projections and replaced
each equation by a pair of opposing half-spaces yielding a consistent
system of hyperslabs. Instead of applying to the system any available
projection method he created around each hyperslab an additional wider
enveloping hyperslab in order to implement a relaxation strategy of
Goffin \cite{goffin1971finite} that advocated interlacing steps of
projection onto the hyperslab's median hyperplane with reflections
into the bounding hyperpalnes of the hyperslabs. Secondly, Bauschke,
Iorio and Koch proposed in \cite{bauschke2014method} an operator,
which they called the intrepid projector, for extending the ART3 method
to handle convex sets in Hilbert space. They replaced the hyperslabs
by enlargements of convex sets. 

We consider the Automatic Relaxation Method (ARM) of \cite{censor1985automatic}
which implemented the strategy of Goffin in a continuous manner, without
using additional enveloping hyperslabs to define whether a point is
``close'' or ``far'' from a hyperslab, as in ART3. The extension
of ART3 from linear hyperslabs to general convex sets by Bauschke,
Iorio and Koch is our inspiration in the present work. We generalize
the ARM algorithm to encompass enlargements of convex sets instead
of being limited to handle only systems of linear hyperslabs.

Enlargements of convex sets are defined as follows.
\begin{definition}
\cite[Definition 2]{bauschke2014method}\textit{\label{def:enlargement}
}Given a nonempty closed convex subset $C$ of a Hilbert space $\mathcal{H}$,
and $\alpha\geq0$, the set 
\begin{equation}
C_{[\alpha]}:=\left\{ x\in\mathcal{H}\mid d_{C}(x)\leq\alpha\right\} 
\end{equation}

\noindent is the \texttt{\textit{$\alpha$-}}enlargement\texttt{ }of
$C$.
\end{definition}
Full details about ART3 and about Bauschke, Iorio and Koch's algorithm
can be found in their papers. For the readers' convenience we give
here only a brief account. In ART3, in addition to the construction
of hyperslabs and enveloping hyperslabs, the ART3 makes three possible
iterative steps: the projection step, the identity step and the reflection
step. The location of a current iterate $x^{k}$ determines the appropriate
step to be taken. If $x^{k}$ is outside the enveloping additional
hyperslab then it is considered to be ``far'' from the original
hyperplane, and, accordingly, a projection step onto the original
hyperplane (which is the median of the hyperslab) will be taken. If
$x^{k}$ is inside the first hyperslab then ART3 will keep it unchanged.
But, if $x^{k}$ is located in the enveloping hyperslab then it is
considered to be ``near'' the original hyperplane and it will be
reflected into the boundary of the hyperslab.
\begin{figure}[H]
\begin{centering}
\includegraphics[scale=0.8]{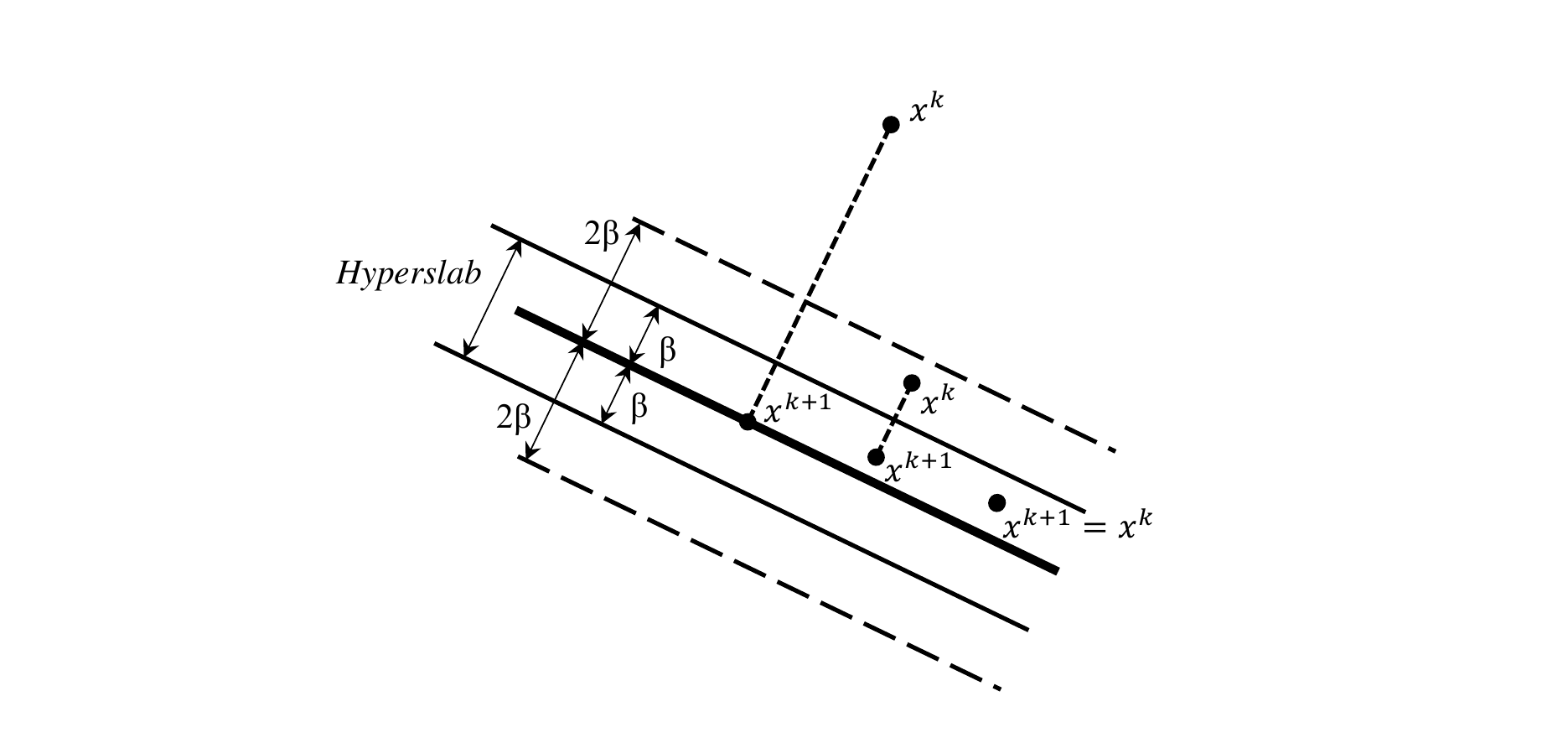}
\par\end{centering}
\caption{The possible iterative steps of ART3. \label{fig:art3}}
\end{figure}

Turning to the work of Bauschke, Iorio and Koch, the operator that
accurately conveys the above ART3 algorithmic strategy to deal with
convex sets is the following operator which we call here the ``geometric
intrepid projector''.
\begin{definition}
\textbf{\textit{\label{def:The-corrected-intrepid}}}\textbf{(The
geometric intrepid projector)} Let $Z$ be a nonempty closed convex
subset of $\mathcal{H}$, let $\beta\in\mathbb{R}$, $\beta>0$ and
set $C:=Z_{[\beta]}$. The geometric intrepid projector $G_{C}:\mathcal{H}\rightarrow\mathcal{H}$
onto $C$ (with respect to $Z$ and $\beta$) is defined by 
\end{definition}
\textit{
\begin{equation}
G_{C}(x):=\begin{cases}
P_{Z}(x), & \textup{if}\,\,\,d_{Z}(x)\geq2\beta,\\
x, & \textup{if}\,\,\,d_{Z}(x)\leq\beta,\\
x+2\left(\dfrac{\beta}{d_{Z}(x)}-1\right)\left(x-P_{Z}(x)\right), & \textup{otherwise.}
\end{cases}\label{eq:8}
\end{equation}
}

The properties of (\ref{eq:8}) for general convex sets remain to
be investigated but we care to present the geometric intrepid projector
$G_{C}$ because Bauschke, Iorio and Koch defined in \cite{bauschke2014method}
a different intrepid projector which does not capture precisely the
ART3 strategy, and which we will name here the ``BIK intrepid projection''.
Contrary to \cite[Example 5]{bauschke2014method}, the third line
in (\ref{eq:BIK intrepid projector}) which is supposed to mimic the
reflection step of ART3, does not actually do so.
\begin{definition}
\textbf{\textit{\label{def:The-original-intrepid} }}\textbf{\cite[Definition 4]{bauschke2014method}
(The Bauschke, Iorio and Koch (BIK) intrepid projector)} Let $Z$
be a nonempty closed convex subset of $\mathcal{H}$, let $\beta\in\mathbb{R}$,
$\beta\geq0$ and set $C:=Z_{[\beta]}$. The projector $Q_{C}:\mathcal{H}\rightarrow\mathcal{H}$
onto $C$ (with respect to $Z$ and $\beta$), henceforth called the
BIK intrepid projector, is defined by {[}for all positive $\beta${]} 
\end{definition}
\noindent \textit{
\begin{equation}
Q_{C}(x):=\begin{cases}
P_{Z}(x), & \textup{if}\,\,\,d_{Z}(x)\geq2\beta,\\
x, & \textup{if}\,\,\,d_{Z}(x)\leq\beta,\\
x+\left(1-\dfrac{d_{Z}(x)}{\beta}\right)\left(x-P_{Z}(x)\right), & \textup{otherwise.}
\end{cases}\label{eq:BIK intrepid projector}
\end{equation}
}

Although it is reasonable to define an algorithm like Bauschke, Iorio
and Koch's algorithm but with $G_{C}$ instead of $Q_{C}$ we do not
stop to do so here and proceed, in the next section, directly to our
valiant projector and its properties.

\section{The Valiant Projector\label{sec:The-valiant-projector}}

In this section we define the valiant projector in Hilbert space and
study its properties.
\begin{definition}
\textbf{(The valiant projector)} Let $Z$ be a nonempty closed convex
subset of $\mathcal{H}$, let $\beta\in\mathbb{R}$, $\beta>0$ and
set $C:=Z_{[\beta]}$. The valiant projector\textit{ }$V_{C}:\mathcal{H}\rightarrow\mathcal{H}$,
onto $C$ (with respect to $Z$ and $\beta$) is defined by 

\noindent 
\begin{equation}
V_{C}(x):=\begin{cases}
x, & \textup{if}\,\,\,d_{Z}(x)\leq\beta,\\
x+\dfrac{\tau}{2}\left(1-\left(\dfrac{\beta}{d_{Z}(x)}\right)^{2}\right)\left(P_{Z}(x)-x\right), & \textup{otherwise,}
\end{cases}\label{eq:the valiant projector}
\end{equation}
with $\tau\in]0,2[$.\footnote{We keep the $\dfrac{\tau}{2}$ in the formula just to show its similarity
with the ARM operator of \cite{censor1985automatic}.}
\end{definition}
The valiant projector works as follows: If the distance of a point
$x$ from the set $Z$ is less than the depth $\beta$ of the enlargement,
i.e., the point is inside the enlargement, then the operator leaves
the point unchanged. Otherwise, if the distance of the point from
the set is greater than the depth of the enlargement, then the valiant
projector brings the point closer to the set $Z$ in the orthogonal
projection direction of the point onto the set. From the term $1-\left(\dfrac{\beta}{d_{Z}(x)}\right)^{2}$
we see that the farther the point is from the set, the operator will
progress towards the set in larger steps. Like in \cite{censor1985automatic},
the overall step-size also incorporates an additional user-chosen
relaxation parameter $\tau$.

The valiant projector has the following useful property.
\begin{proposition}
\label{prop:G is SQNE}Let $Z$ be a nonempty, closed and convex subset
of $\mathcal{H}$, let $\beta\in\mathbb{R}$, $\beta>0.$ If $C:=Z_{[\beta]}$
then the valiant projector $V_{C}$ of (\ref{eq:the valiant projector})
is SQNE.
\end{proposition}
\begin{svmultproof}
We prove that $V_{C}$ is SQNE with respect to $C$ and $\beta$.
For the case $d_{Z}(x)\leq\beta$ we have $V_{C}(x)=x$ and the SQNE
of $V_{C}$ is trivial. For the case $d_{Z}(x)>\beta$, the operator
$V_{C}$ can be written as 
\begin{equation}
V_{C}(x)=\left(1-\frac{\tau}{2}\left(1-\left(\frac{\beta}{d_{Z}(x)}\right)^{2}\right)\right)x+\frac{\tau}{2}\left(1-\left(\frac{\beta}{d_{Z}(x)}\right)^{2}\right)P_{Z}(x),
\end{equation}

\noindent with $\tau\in]0,2[$. Let 
\begin{equation}
\gamma(x):=\frac{\tau}{2}\left(1-\left(\frac{\beta}{d_{Z}(x)}\right)^{2}\right),
\end{equation}

\noindent then we have 
\begin{equation}
V_{C}\left(x\right)=\left(1-\gamma(x)\right)x+\gamma(x)P_{Z}\left(x\right)\label{eq:simple form}
\end{equation}

\noindent with $\gamma(x)\in]0,1[$. Using Proposition \ref{prop:relaxed projector}
with $\varOmega=Z$, $\lambda=\gamma(x)$ and $R=V_{C}$ we have,
for all $c\in Z$,
\begin{equation}
\Vert x-c\Vert^{2}-\Vert V_{C}(x)-c\Vert^{2}\geq\frac{2-\gamma(x)}{\gamma(x)}\Vert x-V_{C}(x)\Vert^{2}\geq\Vert x-V_{C}(x)\Vert^{2}.\label{eq:14}
\end{equation}
Thus, $V_{C}$ is $1$-SQNE therefore, SQNE, by Definition \ref{def:Fixed point-operator}(iv).
\end{svmultproof}

Since both the identity operator $\textup{Id}$ and the projection
$P_{Z}$ are NE, any convex combination of them will be also NE. However,
the dependence of $\gamma$ in (\ref{eq:simple form}) on $x$ requires
special attention when attempting to show that the valiant operator
of (\ref{eq:simple form}) is NE. This is done in the next proposition.
\begin{proposition}
\label{prop:G=00003DNE}Let $Z$ be a nonempty closed convex subset
of $\mathcal{H}$, let $\beta\in\mathbb{R}$, $\beta>0$. If $C:=Z_{[\beta]}$
then the valiant projector $V_{C}$ is NE.
\end{proposition}
\begin{svmultproof}
The proof is split into the three possibilities that have to be considered
according to whether $V_{C}$ realizes the first or the second line
of its definition (\ref{eq:the valiant projector}). 

\textbf{\uline{Possibility A}}: The first line of (\ref{eq:the valiant projector})
holds for both points $x$ and $y$. In this case $d_{Z}(x)\leq\beta$
and $d_{Z}(y)\leq\beta$ thus, $V_{C}(x)=x$ and $V_{C}(y)=y$ so
that $V_{C}$ is trivially NE. 

\textbf{\uline{Possibility B}}: The second line of (\ref{eq:the valiant projector})
holds for both points $x$ and $y$. Without loss of generality, take
any two points $x,y\in\mathcal{H}$ such that 
\begin{equation}
\left\Vert P_{Z}(y)-y\right\Vert \leq\left\Vert P_{Z}(x)-x\right\Vert ,
\end{equation}
and denote 
\begin{align}
d_{Z}(x)=\left\Vert P_{Z}(x)-x\right\Vert =a\beta,\label{eq:16}
\end{align}

\begin{align}
d_{Z}(y)=\left\Vert P_{Z}(y)-y\right\Vert =b\beta,\label{eq:17}
\end{align}

for some real $1<b\leq a.$ Then

\begin{equation}
\frac{\left\Vert V_{C}(x)-x\right\Vert }{\left\Vert V_{C}(y)-y\right\Vert }=\left(\frac{1-\dfrac{1}{a^{2}}}{1-\dfrac{1}{b^{2}}}\right)\frac{\left\Vert P_{Z}(x)-x\right\Vert }{\left\Vert P_{Z}(y)-y\right\Vert },
\end{equation}
so that
\begin{equation}
\left\Vert V_{C}(y)-y\right\Vert \leq\left\Vert V_{C}(x)-x\right\Vert .\label{eq:Ds of Ps}
\end{equation}

There are three possible locations of the points $x,y,V_{C}(x)$ and
$V_{C}(y)$ with respect to the set $Z$, see, Figure \ref{fig:Cases-of-the proof}.
\begin{figure}[H]
\begin{centering}
\includegraphics[scale=0.8]{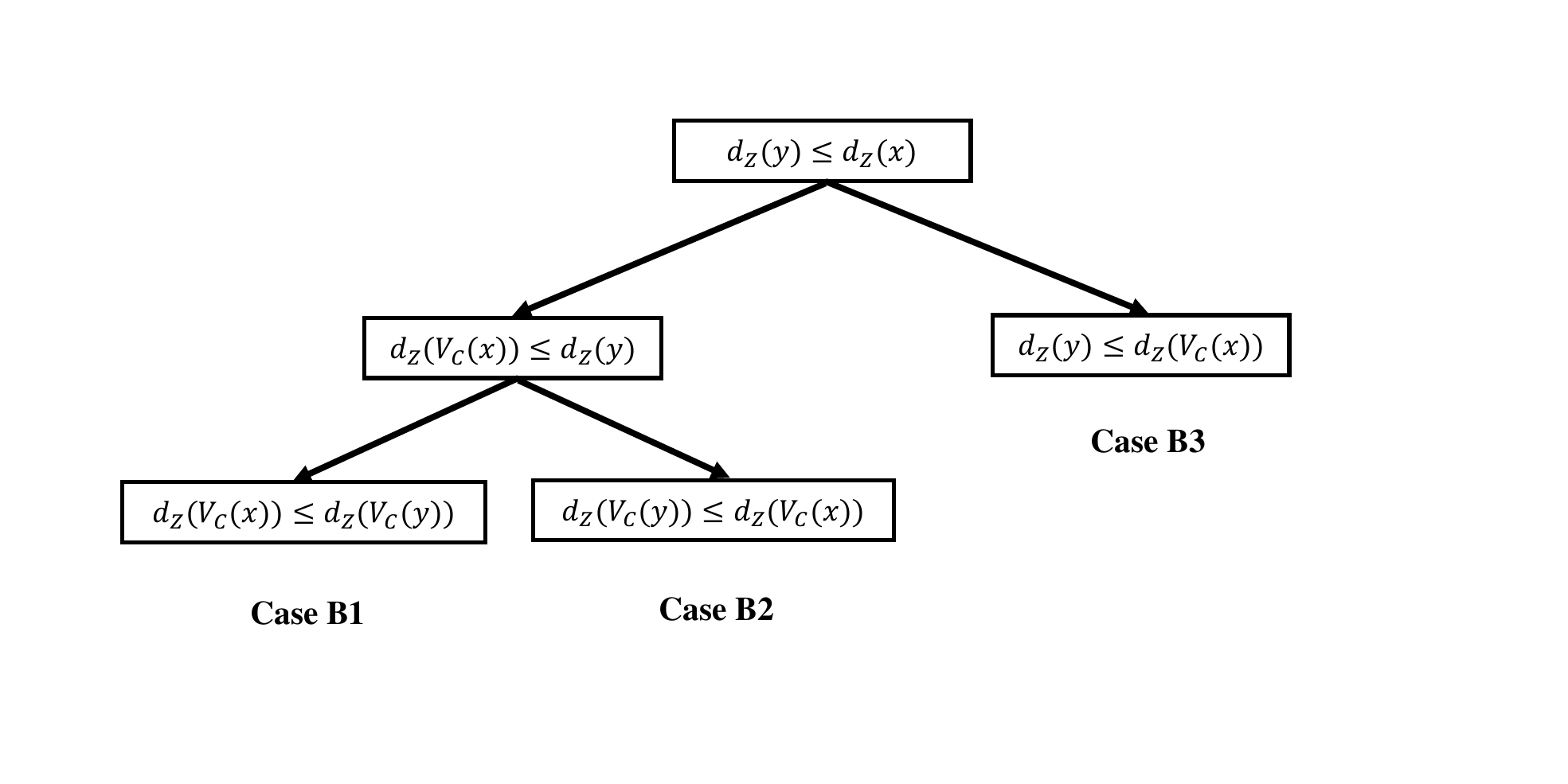}
\par\end{centering}
\caption{Cases for Possibility B in the proof of Proposition \ref{prop:G=00003DNE}
\label{fig:Cases-of-the proof}}
\end{figure}
 Below we discuss each case separately.

\noindent \textbf{Case B1. }Assume that 
\begin{equation}
\left\Vert P_{Z}(V_{C}(x))-V_{C}(x)\right\Vert \leq\left\Vert P_{Z}(V_{C}(y))-V_{C}(y)\right\Vert .
\end{equation}

This implies that

\noindent 
\begin{equation}
\left\Vert P_{Z}(x)-V_{C}(x)\right\Vert \leq\left\Vert P_{Z}(y)-V_{C}(y)\right\Vert ,
\end{equation}

because $P_{Z}(V_{C}(x))=P_{Z}(x)$ and $P_{Z}(V_{C}(y))=P_{Z}(y)$
which follows from the fact that $x$, $V_{C}(x)$ and $P_{Z}(V_{C}(x))$
lie on the same line, and similarly for the other equality, see Figure
\ref{fig:Case-1}.

To study this case we add two enlargements to the set $Z$, one with
a width of $\left\Vert P_{Z}(y)-y\right\Vert $ and the other with
a width of $\left\Vert P_{Z}(y)-V_{C}(y)\right\Vert $. We denote
the intersection point of the line through $x$ and $P_{Z}(x)$ with
the boundary of the first enlargement by $\hat{y}$. Therefore, the
intersection of the above mentioned line with the boundary of the
second enlargement is exactly $V_{C}(\hat{y})$. This is so because,
by (\ref{eq:the valiant projector}), points which are at the same
distance from $Z$ have their images under a valiant operator also
at equal distances from $Z$. See Figure \ref{fig:Case-1}, where
the dashed lines are the enlargements. 
\begin{figure}[H]
\centering{}\includegraphics[scale=0.8]{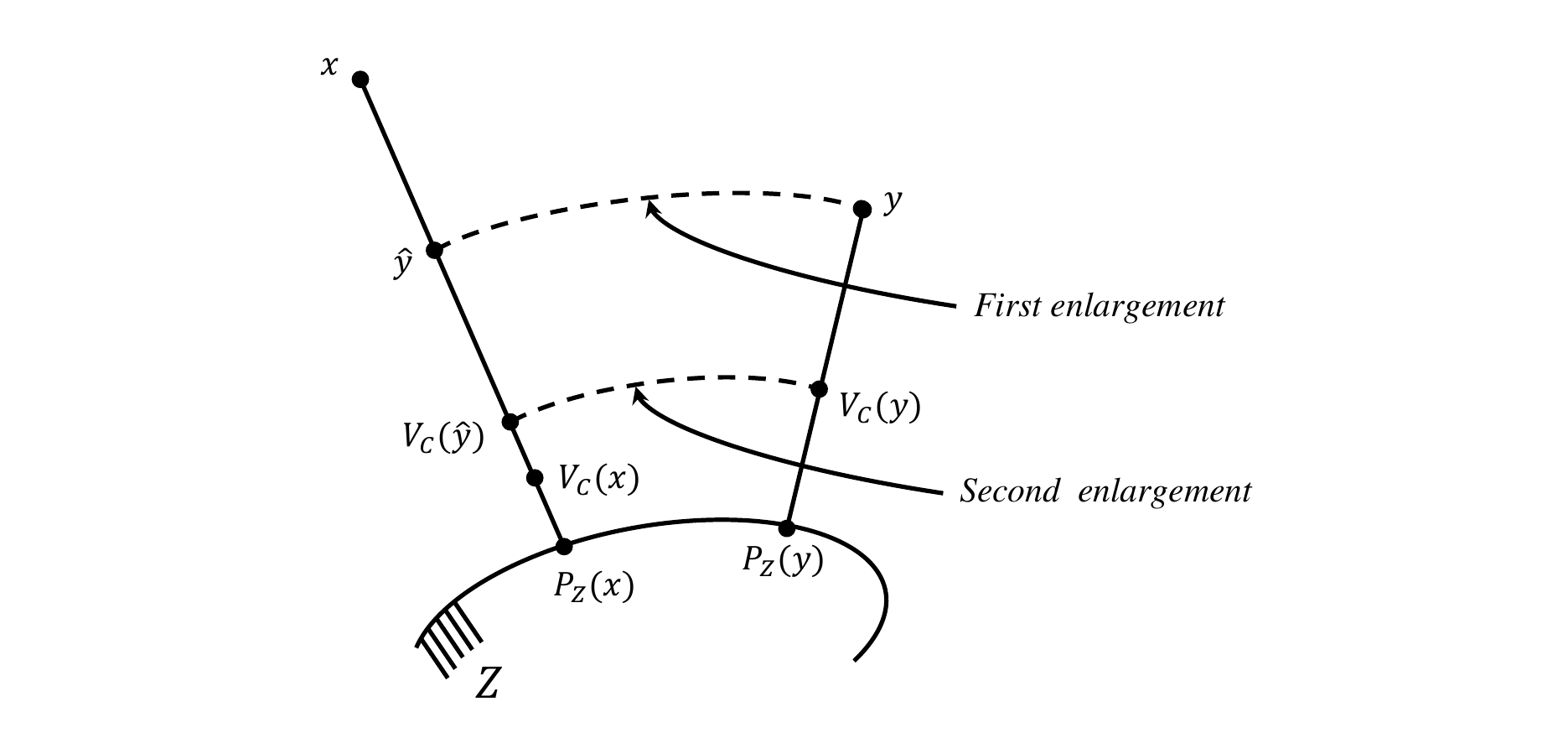}\caption{Case B1 in the proof of Proposition \ref{prop:G=00003DNE}.\label{fig:Case-1}}
\end{figure}

Now we calculate the relevant distances. First, note that
\begin{align}
\left\Vert x-\hat{y}\right\Vert  & =\left\Vert P_{Z}(x)-x\right\Vert -\left\Vert P_{Z}(x)-\hat{y}\right\Vert \nonumber \\
 & =\left\Vert P_{Z}(x)-x\right\Vert -\left\Vert P_{Z}(y)-y\right\Vert \nonumber \\
 & =a\beta-b\beta\nonumber \\
 & =(a-b)\beta.\label{eq:distance 1}
\end{align}
Secondly, 
\begin{align}
\left\Vert V_{C}(x)-V_{C}(\hat{y})\right\Vert = & \left\Vert P_{Z}(y)-V_{C}(y)\right\Vert -\left\Vert P_{Z}(x)-V_{C}(x)\right\Vert \nonumber \\
= & \left\Vert P_{Z}(y)-y\right\Vert -\left\Vert V_{C}(y)-y\right\Vert -\left(\left\Vert P_{Z}(x)-x\right\Vert -\left\Vert V_{C}(x)-x\right\Vert \right)\nonumber \\
= & b\beta-\left|\frac{\tau}{2}\left(1-\left(\frac{\beta}{b\beta}\right)^{2}\right)\right|\left\Vert P_{Z}(y)-y\right\Vert -a\beta+\left|\frac{\tau}{2}\left(1-\left(\frac{\beta}{a\beta}\right)^{2}\right)\right|\left\Vert P_{Z}(x)-x\right\Vert \nonumber \\
= & b\beta-\frac{\tau}{2}\left(1-\frac{1}{b^{2}}\right)b\beta-a\beta+\frac{\tau}{2}\left(1-\frac{1}{a^{2}}\right)a\beta=\beta(a-b)\left(\frac{\tau}{2}\left(1+\frac{1}{ab}\right)-1\right).
\end{align}
Since, by the definitions of the parameters $a,\:b$ and $\tau,$
\begin{equation}
\frac{\tau}{2}\left(1+\frac{1}{ab}\right)-1<1,\label{eq:namely 1}
\end{equation}
we have, by (\ref{eq:distance 1}), that
\begin{equation}
\left\Vert V_{C}(x)-V_{C}(\hat{y})\right\Vert \leq\beta(a-b)=\left\Vert x-\hat{y}\right\Vert .\label{eq:namely 3}
\end{equation}

Since the enlargement of a convex set is also a convex set, 
\begin{equation}
V_{C}(y)=P_{Z_{\left[\left\Vert P_{Z}(y)-V_{C}(y)\right\Vert \right]}}(y),\label{eq:valiant on enlargement 1}
\end{equation}
and 
\begin{equation}
V_{C}(\hat{y})=P_{Z_{\left[\left\Vert P_{Z}(y)-V_{C}(y)\right\Vert \right]}}(\hat{y}),\label{eq:valiant on enlargement 2}
\end{equation}
and 
\begin{equation}
\hat{y}=P_{Z_{\left[\left\Vert P_{Z}(y)-y\right\Vert \right]}}(x).
\end{equation}

As is well-known, the metric projection is NE, so, by (\ref{eq:valiant on enlargement 1})
and (\ref{eq:valiant on enlargement 2}),
\begin{equation}
\left\Vert V_{C}(y)-V_{C}(\hat{y})\right\Vert \leq\left\Vert y-\hat{y}\right\Vert .\label{eq:NE of MP}
\end{equation}

By the characterization of the metric projection, see, e.g., \cite[Theorem 1.2.4]{Ceg-book},
for $P_{Z_{\left[\left\Vert P_{Z}(y)-y\right\Vert \right]}}(x)$ we
have 
\begin{equation}
\left\langle x-\hat{y},y-\hat{y}\right\rangle \leq0,\label{eq:characterization}
\end{equation}
and for $P_{Z_{\left[\left\Vert P_{Z}(y)-V_{C}(y)\right\Vert \right]}}(\hat{y})$
we have 
\begin{equation}
\left\langle x-V_{C}(\hat{y}),V_{C}(y)-V_{C}(\hat{y})\right\rangle \leq0,
\end{equation}
thus, 
\begin{equation}
\left\langle V_{C}(x)-V_{C}(\hat{y}),V_{C}(y)-V_{C}(\hat{y})\right\rangle \geq0.\label{eq:result of characterization}
\end{equation}

We also have,By Proposition \ref{prop:G is SQNE}
\begin{align}
\left\Vert x-y\right\Vert ^{2}= & \left\Vert x-\hat{y}-(y-\hat{y})\right\Vert ^{2}\nonumber \\
= & \left\Vert x-\hat{y}\right\Vert ^{2}+\left\Vert y-\hat{y}\right\Vert ^{2}-2\left\langle x-\hat{y},y-\hat{y}\right\rangle ,\label{eq:inner product 1}
\end{align}
and, similarly, 
\begin{align}
\left\Vert V_{C}(x)-V_{C}(y)\right\Vert ^{2}= & \left\Vert V_{C}(x)-V_{C}(\hat{y})-\left(V_{C}(y)-V_{C}(\hat{y})\right)\right\Vert ^{2}\nonumber \\
= & \left\Vert V_{C}(x)-V_{C}(\hat{y})\right\Vert ^{2}+\left\Vert V_{C}(y)-V_{C}(\hat{y})\right\Vert ^{2}\nonumber \\
 & -2\left\langle V_{C}(x)-V_{C}(\hat{y}),V_{C}(y)-V_{C}(\hat{y})\right\rangle .\label{eq:inner product 2}
\end{align}

Using (\ref{eq:namely 3}), (\ref{eq:NE of MP}), (\ref{eq:characterization})
and (\ref{eq:result of characterization}) in (\ref{eq:inner product 1})
and (\ref{eq:inner product 2}) we get 
\begin{equation}
\left\Vert V_{C}(x)-V_{C}(y)\right\Vert \leq\left\Vert x-y\right\Vert ,\label{eq:V(Z) is NE}
\end{equation}
which proves the nonexpansivity of $V_{C}$ in this case.

\noindent \textbf{Case B2. }With an argument similar to the argument
at the beginning of Case B1 we can assume here that \textit{
\begin{equation}
\left\Vert P_{Z}(y)-V_{C}(y)\right\Vert \leq\left\Vert P_{Z}(x)-V_{C}(x)\right\Vert .
\end{equation}
}

To study this case we add three enlargements to the set $Z$, one
with a width of $\left\Vert P_{Z}(y)-y\right\Vert $, the second with
a width of $\left\Vert P_{Z}(x)-V_{C}(x)\right\Vert $ and the third
with a width of $\left\Vert P_{Z}(y)-V_{C}(y)\right\Vert $. We denote
the intersection point between the line through $x$ and $P_{Z}(x)$
with the boundary of the first enlargement by $\hat{y}$. Therefore,
the intersection of the above mentioned line with the boundary of
the third enlargement is exactly $V_{C}(\hat{y})$. See Figure \ref{fig:Case-2},
where the dashed lines are the enlargements. 
\begin{figure}[H]
\begin{centering}
\includegraphics[scale=0.8]{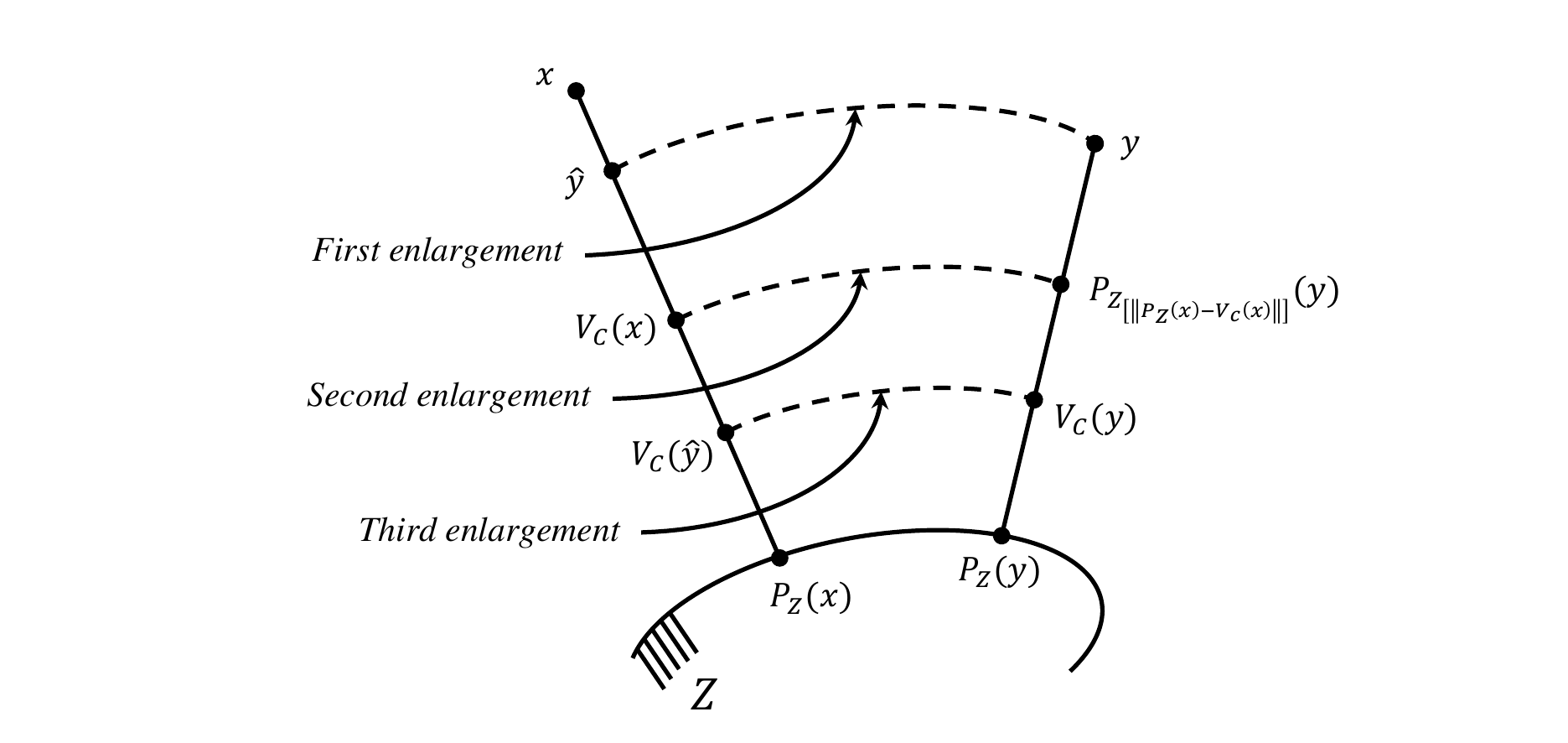}
\par\end{centering}
\caption{Case B2 in the proof of Proposition \ref{prop:G=00003DNE}.\label{fig:Case-2}}
\end{figure}

Now, using a reasoning similar to that in (\ref{eq:namely 1})\textendash (\ref{eq:namely 3}),
we calculate 
\begin{align}
\left\Vert V_{C}(\hat{y})-V_{C}(x)\right\Vert = & \left\Vert P_{Z}(x)-V_{C}(x)\right\Vert -\left\Vert P_{Z}(x)-V_{C}(\hat{y})\right\Vert \nonumber \\
= & \beta(a-b)\left(1-\frac{\tau}{2}\left(1+\frac{1}{ab}\right)\right)\leq\beta(a-b)=\left\Vert x-\hat{y}\right\Vert .\label{eq:as in case 1}
\end{align}

The intersection point between the second enlargement and the line
which passes through $y$ and $P_{Z}(y)$ is $P_{Z_{\left[\left\Vert P_{Z}(x)-V_{C}(x)\right\Vert \right]}}(y)$.
By the nonexpansivity of the metric projection we have 
\begin{equation}
\left\Vert P_{Z_{\left[\left\Vert P_{Z}(x)-V_{C}(x)\right\Vert \right]}}(y)-V_{C}(x)\right\Vert \leq\left\Vert y-\hat{y}\right\Vert .\label{eq:NE of MP2}
\end{equation}

By the characterization of the metric projection, for $P_{Z_{\left[\left\Vert P_{Z}(x)-V_{C}(x)\right\Vert \right]}}(y)$
we have 
\begin{equation}
\left\langle y-P_{Z_{\left[\left\Vert P_{Z}(x)-V_{C}(x)\right\Vert \right]}}(y),V_{C}(x)-P_{Z_{\left[\left\Vert P_{Z}(x)-V_{C}(x)\right\Vert \right]}}(y)\right\rangle \leq0,
\end{equation}
and so 
\begin{equation}
\left\langle V_{C}(y)-P_{Z_{\left[\left\Vert P_{Z}(x)-V_{C}(x)\right\Vert \right]}}(y),V_{C}(x)-P_{Z_{\left[\left\Vert P_{Z}(x)-V_{C}(x)\right\Vert \right]}}(y)\right\rangle \geq0.\label{eq:second result of characterization}
\end{equation}

Now by (\ref{eq:characterization}), (\ref{eq:as in case 1}), (\ref{eq:NE of MP2}),
(\ref{eq:second result of characterization}), by using similar calculations
as in (\ref{eq:inner product 1}), and by replacing $V_{C}(\hat{y})$
by $P_{Z_{\left[\left\Vert P_{Z}(x)-V_{C}(x)\right\Vert \right]}}(y)$
in (\ref{eq:inner product 2}) we obtain (\ref{eq:V(Z) is NE}), namely,
the nonexpansivity of $V_{C}$.

\noindent \textbf{Case B3.} With an argument similar to the argument
at the beginning of Case B1 we can assume here that\textit{
\begin{equation}
\left\Vert P_{Z}(y)-y\right\Vert \leq\left\Vert P_{Z}(x)-V_{C}(x)\right\Vert .
\end{equation}
}

Consult Figure \ref{fig:Case-3}. The points $x,y,P_{Z}(x),P_{Z}(y),V_{C}(x)$
and $V_{C}(y)$ depict the situation for this case. To study this
case we add three enlargements to the set $Z$, one with a width of
$\left\Vert P_{Z}(x)-x\right\Vert $, the second with a width of $\left\Vert P_{Z}(y)-y\right\Vert $
and the third with a width of $\left\Vert P_{Z}(y)-V_{C}(y)\right\Vert $.
We denote the intersection between the line through $x$ and $P_{Z}(x)$
and the boundary of the second enlargement by $\hat{y}$. Therefore,
as argued earlier, the intersection of this line with the boundary
of the third enlargement is precisely $V_{C}(\hat{y})$. We denote
the intersection between the line through $y$ and $P_{Z}(y)$ with
the boundary of the first enlargement by $\hat{x}$. 
\begin{figure}[H]
\begin{centering}
\includegraphics[scale=0.8]{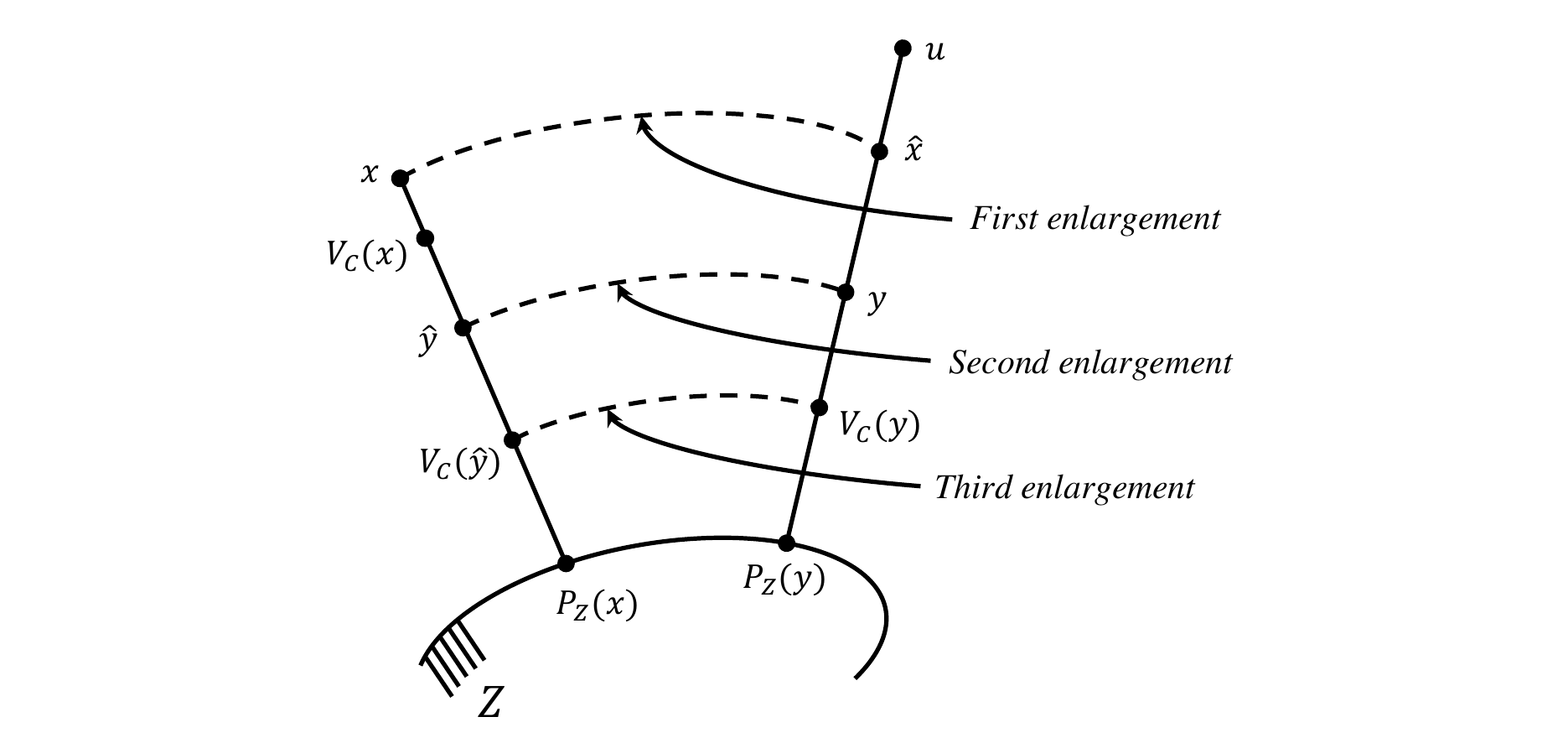}
\par\end{centering}
\centering{}\caption{Case B3 in the proof of Proposition \ref{prop:G=00003DNE}.\label{fig:Case-3}}
\end{figure}

Now we calculate
\begin{align}
\left\Vert x-V_{C}(y)\right\Vert ^{2}= & \left\Vert x-V_{C}(\hat{y})-\left(V_{C}(y)-V_{C}(\hat{y})\right)\right\Vert ^{2}\nonumber \\
= & \left\Vert x-V_{C}(\hat{y})\right\Vert ^{2}+\left\Vert V_{C}(y)-V_{C}(\hat{y})\right\Vert ^{2}\nonumber \\
 & -2\left\langle x-V_{C}(\hat{y}),V_{C}(y)-V_{C}(\hat{y})\right\rangle .\label{eq:InPr-1}
\end{align}
Calculating the left-hand side of (\ref{eq:InPr-1}) a bit differently
we may write
\begin{align}
\left\Vert x-V_{C}(y)\right\Vert ^{2}= & \left\Vert x-\hat{x}-\left(V_{C}(y)-\hat{x}\right)\right\Vert ^{2}\nonumber \\
= & \left\Vert x-\hat{x}\right\Vert ^{2}+\left\Vert V_{C}(y)-\hat{x}\right\Vert ^{2}-2\left\langle x-\hat{x},V_{C}(y)-\hat{x}\right\rangle .\label{eq:InPr-3}
\end{align}
Finally, we write
\begin{align}
\left\Vert x-y\right\Vert ^{2}= & \left\Vert x-\hat{x}-(y-\hat{x})\right\Vert ^{2}\nonumber \\
= & \left\Vert x-\hat{x}\right\Vert ^{2}+\left\Vert y-\hat{x}\right\Vert ^{2}-2\left\langle x-\hat{x},y-\hat{x}\right\rangle .\label{eq:InPr-4}
\end{align}
By subtracting (\ref{eq:inner product 2}) from (\ref{eq:InPr-1}),
by subtracting (\ref{eq:InPr-4}) from (\ref{eq:InPr-3}), and then
by subtracting the second result from the first and because 
\begin{equation}
\left\Vert x-V_{C}(\hat{y})\right\Vert ^{2}=\left\Vert V_{C}(y)-\hat{x}\right\Vert ^{2},
\end{equation}
 we obtain,
\begin{align}
\left\Vert x-y\right\Vert ^{2}-\left\Vert V_{C}(x)-V_{C}(y)\right\Vert ^{2}= & \left\Vert y-\hat{x}\right\Vert ^{2}-\left\Vert V_{C}(x)-V_{C}(\hat{y})\right\Vert ^{2}\nonumber \\
 & +2\left\langle V_{C}(x)-V_{C}(\hat{y}),V_{C}(y)-V_{C}(\hat{y})\right\rangle \nonumber \\
 & -2\left\langle x-V_{C}(\hat{y}),V_{C}(y)-V_{C}(\hat{y})\right\rangle \nonumber \\
 & +2\left\langle x-\hat{x},V_{C}(y)-\hat{x}\right\rangle -2\left\langle x-\hat{x},y-\hat{x}\right\rangle .\label{eq:Result-Case3}
\end{align}
Along the line through $x$ and $P_{Z}(x)$ we have, by using (\ref{eq:Ds of Ps}),
\begin{align}
\left\Vert V_{C}(\hat{y})-V_{C}(x)\right\Vert = & \left\Vert V_{C}(\hat{y})-\hat{y}\right\Vert +\left\Vert \hat{y}-V_{C}(x)\right\Vert \nonumber \\
= & \left\Vert V_{C}(y)-y\right\Vert +\left\Vert \hat{y}-V_{C}(x)\right\Vert \nonumber \\
\leq & \left\Vert V_{C}(x)-x\right\Vert +\left\Vert \hat{y}-V_{C}(x)\right\Vert \nonumber \\
= & \left\Vert \hat{y}-x\right\Vert \nonumber \\
= & \left\Vert y-\hat{x}\right\Vert .\label{eq:F1}
\end{align}
By the linearity of the inner product and by using the characterization
of the metric projection we have 
\begin{align}
\left\langle x-V_{C}(\hat{y}),V_{C}(y)-V_{C}(\hat{y})\right\rangle \nonumber \\
= & \left\langle x-V_{C}(x),V_{C}(y)-V_{C}(\hat{y})\right\rangle \nonumber \\
 & +\left\langle V_{C}(x)-V_{C}(\hat{y}),V_{C}(y)-V_{C}(\hat{y})\right\rangle \nonumber \\
\leq & \left\langle V_{C}(x)-V_{C}(\hat{y}),V_{C}(y)-V_{C}(\hat{y})\right\rangle .\label{eq:F2}
\end{align}
Finally, let $u$ be a point on the line through $y$ and $P_{Z}(y)$
such that 
\begin{equation}
P_{Z_{\left[\left\Vert P_{Z}(y)-\hat{x}\right\Vert \right]}}(u)=\hat{x.}
\end{equation}
By the characterization of the metric projection we have 
\begin{equation}
\left\langle x-\hat{x},u-\hat{x}\right\rangle \leq0,
\end{equation}
so, 
\begin{equation}
\left\langle y-\hat{x},x-\hat{x}\right\rangle \geq0.\label{eq:MP3}
\end{equation}
By the linearity of the inner product and by (\ref{eq:MP3}) we have
\begin{align}
\left\langle V_{C}(y)-\hat{x},x-\hat{x}\right\rangle = & \left\langle V_{C}(y)-y,x-\hat{x}\right\rangle +\left\langle y-\hat{x},x-\hat{x}\right\rangle \nonumber \\
\geq & \left\langle y-\hat{x},x-\hat{x}\right\rangle .\label{eq:F3}
\end{align}
Using (\ref{eq:F1}), (\ref{eq:F2}) and (\ref{eq:F3}) in (\ref{eq:Result-Case3}),
we have
\begin{equation}
\left\Vert V_{C}(x)-V_{C}(y)\right\Vert \leq\left\Vert x-y\right\Vert .
\end{equation}
 By (\ref{eq:Ds of Ps}) and considerations as in Cases B1 and B2
we get the nonexpansivity of $V_{C}$, and the proof is complete.

\textbf{\uline{Possibility C}}\textbf{: }The first line of (\ref{eq:the valiant projector})
holds for $y$, i.e., $V_{C}(y)=y$, and the second line of (\ref{eq:the valiant projector})
holds for $x.$

Now we discuss the situation 
\begin{align}
d_{Z}(x)>\beta,\\
d_{Z}(y)\leq\beta,
\end{align}
and repeat (\ref{eq:16}) and (\ref{eq:17}) but this time for some
real $0<b\leq1<a.$

\textbf{Case C1.} With an argument similar to the argument at the
beginning of Case B1 we can assume here that\textit{
\begin{equation}
\left\Vert P_{Z}(x)-V_{C}(x)\right\Vert \geq\left\Vert P_{Z}(y)-y\right\Vert .
\end{equation}
}To study this case we add two enlargements to the set $Z$, one with
a width of $\left\Vert P_{Z}(y)-y\right\Vert $ and the other with
a width of $\left\Vert P_{Z}(x)-V_{C}(x)\right\Vert $. We denote
the intersection point of the line through $x$ and $P_{Z}(x)$ with
the boundary of the first enlargement by $\hat{y}$, and with the
boundary of the second enlargement by $V_{C}(x)$. We also denote
the intersection point of the line through $y$ and $P_{Z}(y)$ with
the boundary of the second enlargement by $V_{C}(\hat{x})$.

\begin{figure}[H]
\begin{centering}
\includegraphics[scale=0.8]{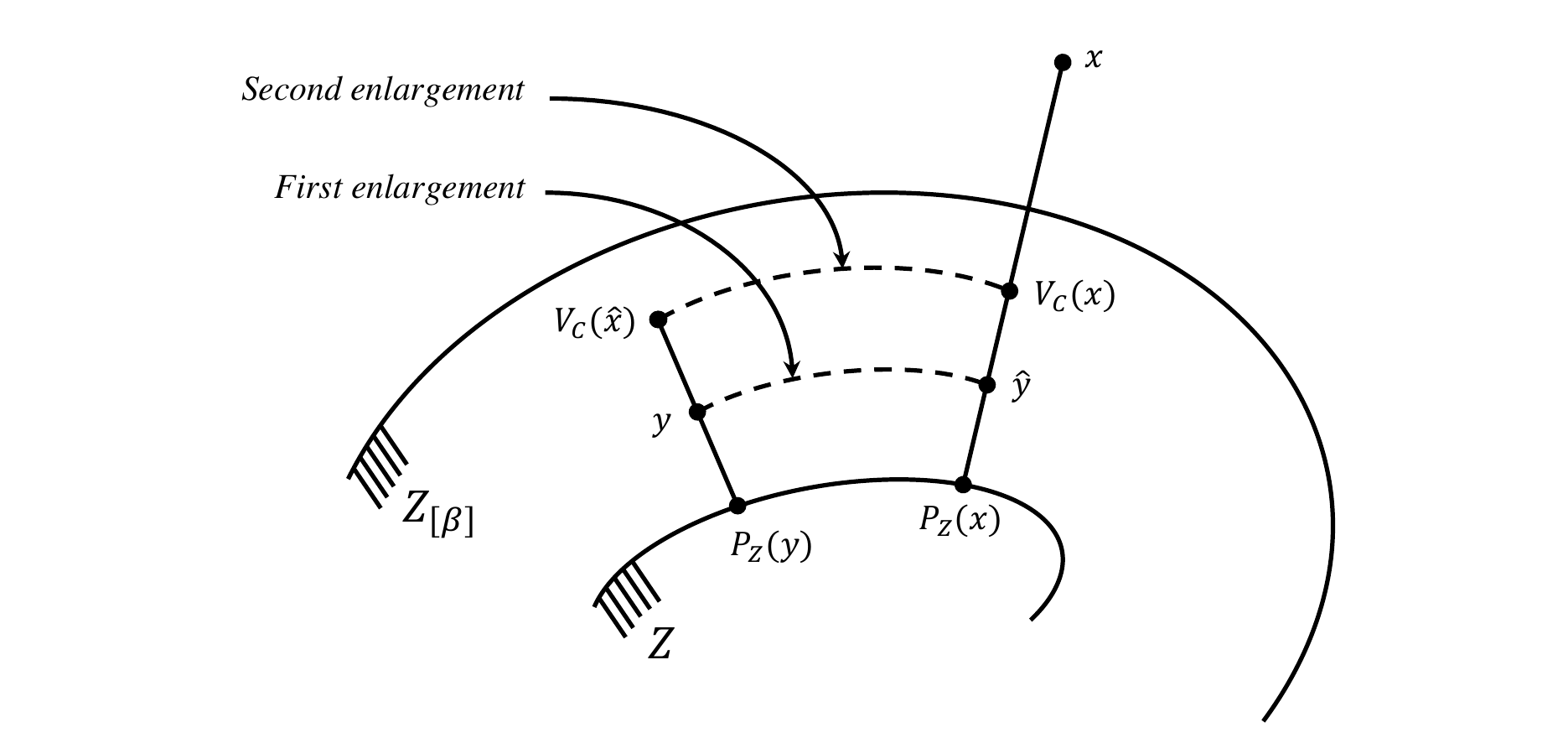}
\par\end{centering}
\centering{}\caption{Case C1 in the proof of Proposition \ref{prop:G=00003DNE}.\label{fig:Case-C1}}
\end{figure}

Now we have \textit{
\begin{equation}
\left\Vert x-\hat{y}\right\Vert \geq\left\Vert V_{C}(x)-\hat{y}\right\Vert .\label{eq:59}
\end{equation}
}We also have 
\begin{align}
\left\Vert V_{C}(x)-V_{C}(y)\right\Vert ^{2}= & \left\Vert V_{C}(x)-y\right\Vert ^{2}\nonumber \\
= & \left\Vert V_{C}(x)-\hat{y}-\left(y-\hat{y}\right)\right\Vert ^{2}\nonumber \\
= & \left\Vert V_{C}(x)-\hat{y}\right\Vert ^{2}+\left\Vert y-\hat{y}\right\Vert ^{2}\nonumber \\
 & -2\left\langle V_{C}(x)-\hat{y},y-\hat{y}\right\rangle .\label{eq:60}
\end{align}
The characterization of the metric projection for $\hat{y}=P_{Z_{\left[\left\Vert P_{Z}(y)-y\right\Vert \right]}}(x)$
allows us to reuse (\ref{eq:characterization}) and also yields
\begin{equation}
\left\langle V_{C}(x)-\hat{y},y-\hat{y}\right\rangle \leq0.\label{eq:61}
\end{equation}
Using (\ref{eq:59}), (\ref{eq:60}), (\ref{eq:characterization}),
(\ref{eq:61}) and (\ref{eq:inner product 1}) proves the nonexpansivity
of $V_{C}$ in this case.

\textbf{Case C2.} With an argument similar to the argument at the
beginning of Case B1 we can assume here that\textit{
\begin{equation}
\left\Vert P_{Z}(x)-V_{C}(x)\right\Vert <\left\Vert P_{Z}(y)-y\right\Vert .
\end{equation}
}
\begin{figure}[H]
\begin{centering}
\includegraphics[scale=0.8]{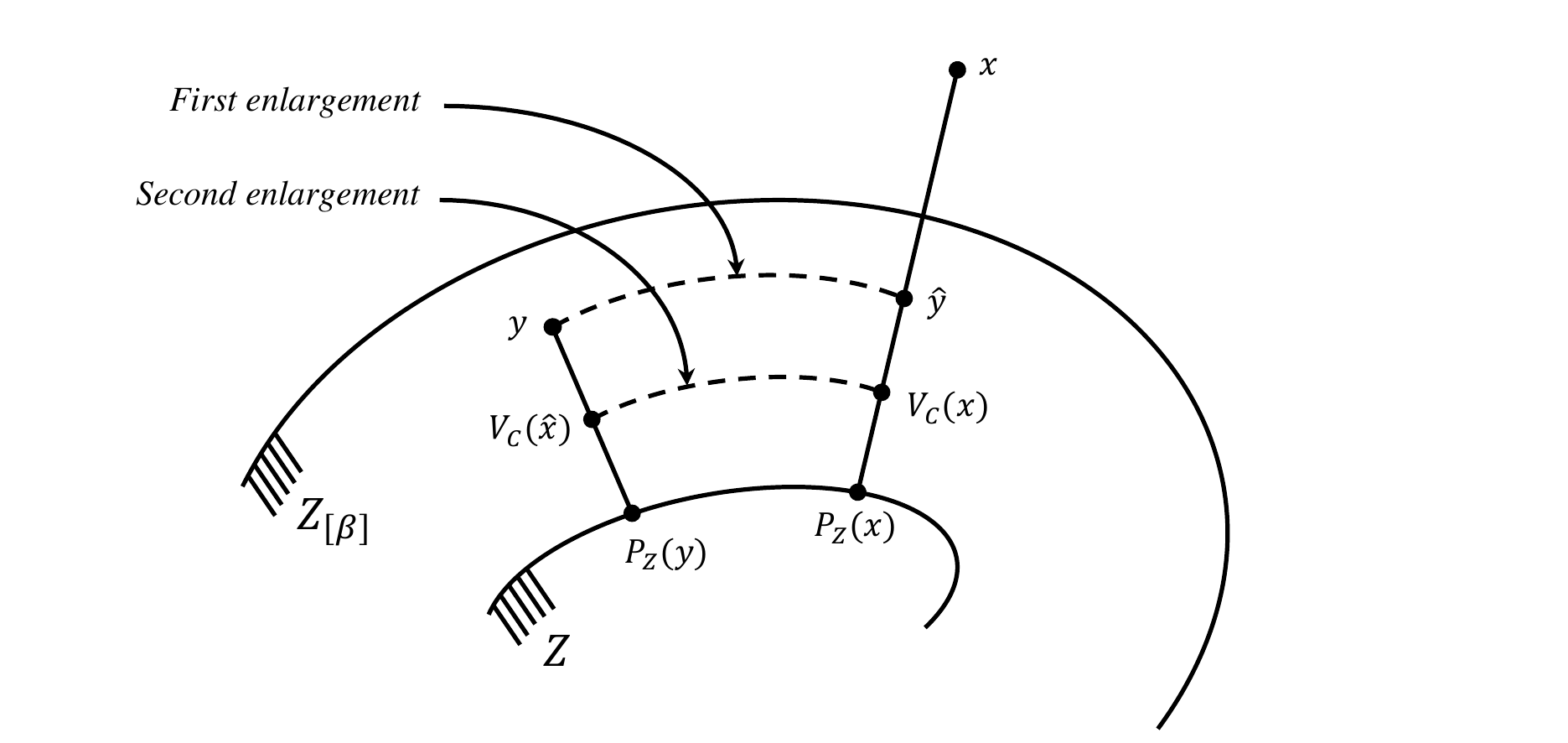}
\par\end{centering}
\caption{Case C2 in the proof of Proposition \ref{prop:G=00003DNE}.\label{fig:Case-C2}}
\end{figure}
To study this case we repeat the addition of two enlargements exactly
as described for Case C1 above. The only difference will be (see Figure
\ref{fig:Case-C2}) that the two enlargements have changed their positions.

The calculation of (\ref{eq:distance 1}) in Case B1 holds again here
verbatim. Now the following holds (recall that $0<b\leq1<a.)$
\begin{align}
\left\Vert V_{C}(x)-\hat{y}\right\Vert = & \left\Vert P_{Z}(y)-y\right\Vert -\left\Vert P_{Z}(y)-V_{C}(\hat{x})\right\Vert \nonumber \\
= & \left\Vert P_{Z}(y)-y\right\Vert -\left(\left\Vert P_{Z}(x)-x\right\Vert -\left\Vert V_{C}(x)-x\right\Vert \right)\nonumber \\
= & b\beta-\left(a\beta-\left|\frac{\tau}{2}\left(1-\left(\frac{\beta}{a\beta}\right)^{2}\right)\right|\left\Vert P_{Z}(x)-x\right\Vert \right)\nonumber \\
= & \frac{\tau}{2}\left(1-\dfrac{1}{a^{2}}\right)a\beta-\left(a-b\right)\beta.\label{eq:63}
\end{align}
By subtracting (\ref{eq:63}) from (\ref{eq:distance 1}) we obtain
\begin{equation}
\left\Vert x-\hat{y}\right\Vert -\left\Vert V_{C}(x)-\hat{y}\right\Vert =2\left(a-b\right)\beta-\frac{\tau}{2}\left(1-\dfrac{1}{a^{2}}\right)a\beta.\label{eq:64}
\end{equation}
Now we calculate the ratio 
\begin{equation}
\dfrac{2\left(a-b\right)\beta}{\dfrac{\tau}{2}\left(1-\dfrac{1}{a^{2}}\right)a\beta}=\dfrac{2}{\tau}\dfrac{a-b}{a-1}2\dfrac{a}{a+1}.\label{eq:65}
\end{equation}
Since $0<\tau<2$ and $0<b\leq1$ we have $\dfrac{2}{\tau}>1$ and
$\dfrac{a-b}{a-1}\geq1$ and we have $\dfrac{1}{2}<\dfrac{a}{a+1}<1,$
so that $1<2\dfrac{a}{a+1}<2,$ and then the ratio in (\ref{eq:65})
is greater than 1, so the right side of (\ref{eq:64}) is greater
than 0 and hence 
\begin{equation}
\left\Vert x-\hat{y}\right\Vert >\left\Vert V_{C}(x)-\hat{y}\right\Vert .\label{eq:66}
\end{equation}
Reusing (\ref{eq:characterization}) we have 
\begin{equation}
\left\langle V_{C}(x)-\hat{y},y-\hat{y}\right\rangle \geq0.\label{eq:67}
\end{equation}
Now by (\ref{eq:inner product 1}), (\ref{eq:characterization}),
(\ref{eq:60}), (\ref{eq:66}) and (\ref{eq:67}) proves the nonexpansivity
of $V_{C}$ in this case.
\end{svmultproof}

Another important feature of the valiant operator is the following.
\begin{proposition}
\label{prop:FixG=00003DC}Let $Z$ be a nonempty closed convex subset
of $\mathcal{H}$, let $\beta\in\mathbb{R}$, $\beta>0.$ If $C:=Z_{[\beta]}$
then the valiant projector $V_{C}$ has the property 
\begin{equation}
\textnormal{Fix}V_{Z_{[\beta]}}=Z_{[\beta]}.
\end{equation}
\end{proposition}
\begin{svmultproof}
If $x\in C$ then $d_{Z}(x)\leq\beta$ and so $V_{C}(x)=x$. Therefore,
$x\in\textnormal{Fix}V_{C}$. If $x\in\textnormal{Fix}V_{C}$ then
$V_{C}(x)=x$. If the first line of (\ref{eq:the valiant projector})
holds then $x\in C$. If the second line of (\ref{eq:the valiant projector})
holds then 
\begin{equation}
x=x+\frac{\tau}{2}\left(1-\left(\frac{\beta}{d_{Z}(x)}\right)^{2}\right)\left(P_{Z}(x)-x\right),
\end{equation}
which leads to either $\beta=d_{Z}(x)$ or $P_{Z}(x)=x$, implying
in both cases that $x\in C$. Therefore, if $x\in\textnormal{Fix}V_{C}(x)$
then $x\in C$.
\end{svmultproof}

\section{The Valiant Projections Method\label{sec:The-method}}

Now we are ready to present our algorithm that employs valiant projections
and prove its convergence. A sequence $\left\{ i_{k}\right\} _{k=0}^{\infty}$
of indices is called a cyclic control sequence on $I:=\left\{ 1,2,\ldots,m\right\} $
if $i_{k}=k\left(\textup{mod}\,m\right)+1$.

Let $Z_{1},Z_{2},\ldots,Z_{m}\subseteq\mathcal{H}$ be nonempty closed
convex sets and $C_{i}=(Z_{i})_{[\beta_{i}]}$ their enlargements
with $\beta_{i}>0$ for all $i\in I$. Assume that $\bigcap_{i\in I}C_{i}\neq\textrm{Ø}$.

\begin{algorithm}[H]
\caption{The Valiant Projection Method (VPM)\textbf{\textit{\label{alg:vpm}}}}

\noindent \textbf{Initialization}: \textit{$x^{0}\in\mathcal{H}$
}is arbitrary.

\noindent \textbf{Iterative Step}: Given the current iterate $x^{k}$,
calculate 
\begin{align}
x^{k+1} & =V_{C_{i_{k}}}(x^{k})\nonumber \\
 & =\begin{cases}
x^{k}, & \textup{if}\,\,\,d_{Z_{i_{k}}}(x)\leq\beta_{i_{k}},\\
x^{k}+\tau_{k}\left(1-\left(\dfrac{\beta_{i_{k}}}{d_{Z_{i_{k}}}(x^{k})}\right)^{2}\right)\left(P_{Z_{i_{k}}}(x^{k})-x^{k}\right), & \textup{otherwise.}
\end{cases}
\end{align}
where $\left\{ i_{k}\right\} _{k=0}^{\infty}$ is cyclic on $I$,
$\beta_{i_{k}}>0$ and $\tau_{k}\in]0,1[$ for all $k\geq0$.
\end{algorithm}

\begin{theorem}
\label{thm:convergence of CRM}Let $Z_{1},Z_{2},\ldots,Z_{m}\subseteq\mathcal{H}$
be nonempty closed convex sets and $C_{i}=(Z_{i})_{[\beta_{i}]}$
their enlargements with $\beta_{i}>0$, for all $i\in I$. Assume
that $\bigcap_{i\in I}C_{i}\neq\textrm{Ø}$. Any sequence $\left\{ x^{k}\right\} _{k=0}^{\infty}$,
generated by Algorithm \ref{alg:vpm}, converges weakly to a point
$x^{*}\in\bigcap_{i\in I}C_{i}$.
\end{theorem}
\begin{svmultproof}
We wish to apply Theorem \ref{thm:3.6.2} and to this end we show
that all assumptions of that theorem hold here. Let $X\subseteq\mathcal{H}$
be a nonempty closed convex subset, we define an operator $S:X\rightarrow X$
by 
\begin{equation}
S:=\prod_{i=1}^{m}V_{C_{i}}=V_{C_{m}}V_{C_{m-1}}\cdots V_{C_{1}}.
\end{equation}
By Proposition \ref{prop:G is SQNE}, each $V_{C_{i}}$ is SQNE, so,
by Proposition \ref{prop:sQNE}(iii), it is sQNE. By Proposition \ref{prop:FixG=00003DC}
we have 
\begin{equation}
\bigcap_{i=1}^{m}\textnormal{Fix}V_{C_{i}}=\bigcap_{i=1}^{m}C_{i}\neq\textrm{Ø}.\label{eq:Fix(comp)=00003DInter(C)}
\end{equation}
 Using Theorem \ref{thm:Fix=00003DintersectionFix} (observe that
this theorem dictates the use of the cyclic control in our algorithm)
with $U_{i}$ as the valiant operators $V_{C_{i}}$, and applying
Proposition \ref{prop:sQNE}(i), we get 
\begin{equation}
\textup{Fix }S=\bigcap_{i=1}^{m}\textnormal{Fix}V_{C_{i}}.\label{eq:FixU=00003DFix(n)}
\end{equation}
Thus, by Remark \ref{rem:SQNE=00003Dcomp+comb}, the operator $S$
is SQNE and so, by Theorem \ref{thm:SQNE=00003DAR}, it is asymptotically
regular. From Proposition \ref{prop:G=00003DNE}, $V_{C_{i}}$ is
NE, and by, Remark \ref{rem:FNE=00003DNE}(ii), $S$ is NE and has
a fixed point according to (\ref{eq:Fix(comp)=00003DInter(C)}) and
(\ref{eq:FixU=00003DFix(n)}). Using the demiclosedness principle
embodied in Theorem \ref{thm:demiclosedness principle} and Definition
\ref{def:demi-closed} for the operator $S$, the operator $S-\textup{Id}$
is demiclosed at 0. 

Since the iterative process of the algorithm consists of repeated
applications of the valiant operator we show next that$\left\{ V_{C_{i_{k}}}\right\} _{k=0}^{\infty}$
is an asymptotically regular sequence of QNE operators. Since $V_{C_{i}}$
is SQNE we have, by Definition \ref{def:Fixed point-operator}(iv)
with $\alpha=1$,
\begin{equation}
\left\Vert V_{C_{i_{k}}}\left(x^{k}\right)-z\right\Vert ^{2}\leq\left\Vert x^{k}-z\right\Vert ^{2}-\left\Vert V_{C_{i_{k}}}\left(x^{k}\right)-x^{k}\right\Vert ^{2},\textup{\;for\;every}\;z\in\textnormal{Fix}V_{C_{i_{k}}},
\end{equation}
which guarantees that
\begin{equation}
\left\Vert x^{k+1}-z\right\Vert ^{2}\leq\left\Vert x^{k}-z\right\Vert ^{2}-\left\Vert x^{k+1}-x^{k}\right\Vert ^{2},\;\textup{for\;every}\;z\in\bigcap_{k=0}^{\infty}\textnormal{Fix}V_{C_{i_{k}}}.
\end{equation}
Consequently, $\left\{ x^{k}\right\} $ is Fejér-monotone with respect
to $\bigcap_{k=0}^{\infty}\textnormal{Fix}V_{C_{i_{k}}}$ thus it
is bounded. Therefore, $\left\{ \left\Vert x^{k}-z\right\Vert \right\} _{k=0}^{\infty}$
is monotonically decreasing thus convergent, which yields 
\begin{equation}
\lim_{k\rightarrow\infty}\left\Vert x^{k+1}-x^{k}\right\Vert =0.\label{eq:ARsequence}
\end{equation}
According to Definition \ref{def:asymp-reg-seq} $\left\{ V_{C_{i_{k}}}\right\} _{k=0}^{\infty}$
is an asymptotically regular sequence of QNE operators. Finally, to
justify (\ref{eq:6}) we compute the following limit using the triangle
inequality and (\ref{eq:ARsequence}). 
\begin{align}
 & \lim_{k}\left\Vert S\left(x^{k}\right)-x^{k}\right\Vert \nonumber \\
 & =\lim_{k}\left\Vert \prod_{i=1}^{m}V_{C_{i}}\left(x^{k}\right)-x^{k}\right\Vert \nonumber \\
 & =\lim_{k}\left\Vert x^{k+m}-x^{k}\right\Vert \nonumber \\
 & =\lim_{k}\left\Vert x^{k+m}-x^{k-1+m}+x^{k-1+m}-x^{k-2+m}+x^{k-2+m}\cdots-x^{k+1}+x^{k+1}-x^{k}\right\Vert \nonumber \\
 & \leq\lim_{k}\left(\left\Vert x^{k+m}-x^{k-1+m}\right\Vert +\left\Vert x^{k-1+m}-x^{k-2+m}\right\Vert +\cdots+\left\Vert x^{k+1}-x^{k}\right\Vert \right)\nonumber \\
 & =\lim_{k}\left\Vert x^{k+m}-x^{k-1+m}\right\Vert +\lim_{k}\left\Vert x^{k-1+m}-x^{k-2+m}\right\Vert +\cdots+\lim_{k}\left\Vert x^{k+1}-x^{k}\right\Vert \nonumber \\
 & =0.
\end{align}
We have proved that all the assumptions of Theorem \ref{thm:3.6.2}
are satisfied. Therefore, $\left\{ x^{k}\right\} _{k=0}^{\infty}$
converges weakly to a point $x^{*}\in\textnormal{Fix}S$ and, by (\ref{eq:Fix(comp)=00003DInter(C)})
and (\ref{eq:FixU=00003DFix(n)}), $x^{*}\in\bigcap_{i=1}^{m}C_{i}$.
\end{svmultproof}

\section{Conclusions\label{sec:Conclusions}}

In Table \ref{tab:Approaches-to-Goffin's} we depict features and
relationships between the algorithmic operators which grew out from
Goffin's principle. The idea of enlargements led to the extensions
of the algorithms ART3 and ARM to handle convex sets. As stated in
\cite[Theorems 11 and 14]{bauschke2014method}, convergence of the
method of cyclic intrepid projections (CycIP), see, \cite[Algorithm 9]{bauschke2014method}
is guaranteed if the interior of the intersection of the sets is not
empty. In the present work with valiant operators this condition is
not required for the convergence of our VPM algorithm.
\begin{table}[H]
\caption{{\small{}Approaches to Goffin's principle\label{tab:Approaches-to-Goffin's}}}

\centering{}\includegraphics[scale=0.8]{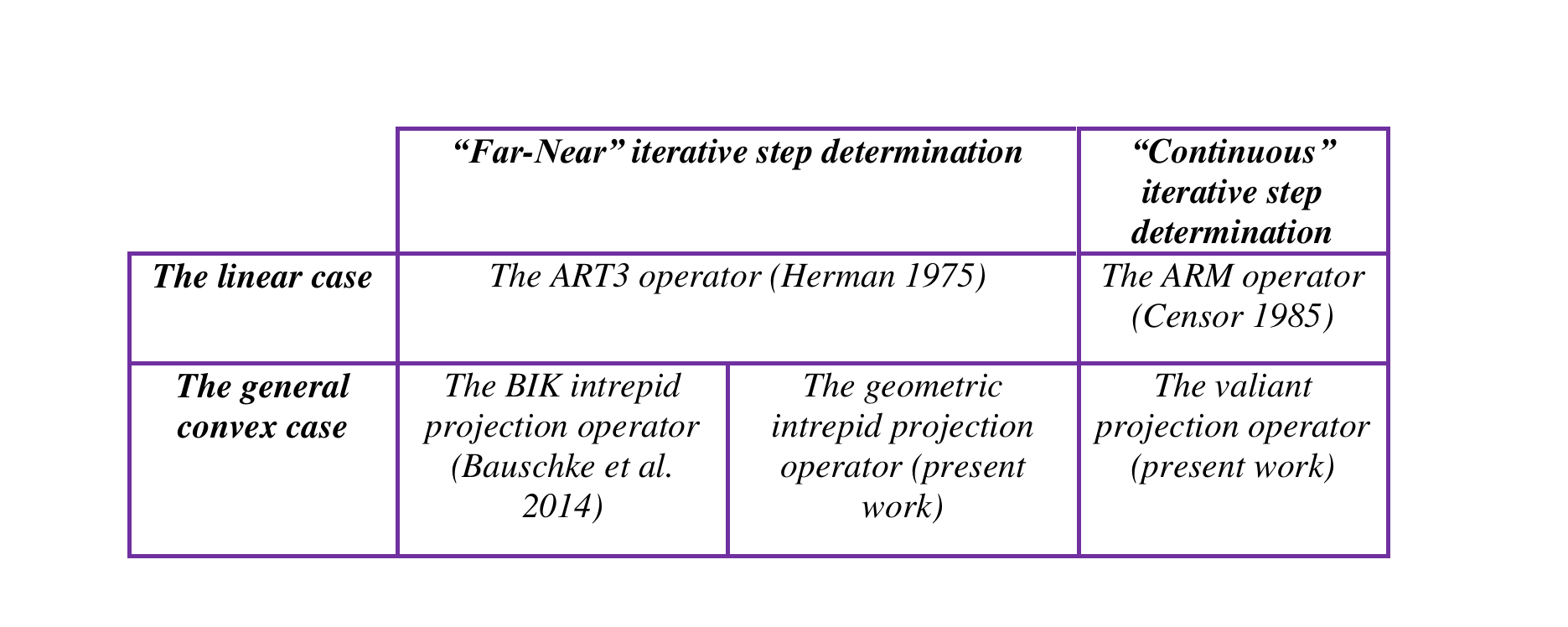}
\end{table}

\begin{acknowledgements}
We thank Tommy Elfving for reading several parts of this paper and
making enlightening comments. We are indebted to the reviewers and
to the Editor-in-Chief Franco Giannessi for their insightful and constructive
comments that helped us improve the paper. This work was supported
by Research Grant No. 2013003 of the United States-Israel Binational
Science Foundation (BSF).
\end{acknowledgements}

\bibliographystyle{unsrt}
\addcontentsline{toc}{section}{\refname}\bibliography{Bib-rafiq-PhD-2016}

\end{document}